\documentclass{svmult}
\usepackage[utf8]{inputenc}
\usepackage{stmaryrd,mathrsfs,bm,amsthm,mathtools,yfonts,amssymb}
\usepackage{color}
\usepackage{makeidx}

\newtheorem{Theorem}{Theorem}[section]
\newtheorem{Proposition}{Proposition}[section]
\newtheorem{Lemma}{Lemma}[section]
\newtheorem{Corollary}{Corollary}[section]
\newtheorem{Definition}{Definition}[section]

\newcommand{\bTheorem}[1]{
\begin{Theorem} \label{T#1} }
\newcommand{\eT}{\end{Theorem}}

\newcommand{\bProposition}[1]{
\begin{Proposition} \label{P#1}}
\newcommand{\eP}{\end{Proposition}}

\newcommand{\bLemma}[1]{
\begin{Lemma} \label{L#1} }
\newcommand{\eL}{\end{Lemma}}

\newcommand{\bCorollary}[1]{
\begin{Corollary} \label{C#1} }
\newcommand{\eC}{\end{Corollary}}

\newcommand{\bFormula}[1]{
\begin{equation} \label{#1}}
\newcommand{\eF}{\end{equation}}

\newcommand{\pat}{\partial_t}
\newcommand{\diver}{\operatorname{div}}

\newcommand{\vr}{\varrho}

\newcommand{\vu}{\vc{u}}
\newcommand{\vc}[1]{{\bf #1}}

\newcommand{\dx}{\,{\rm d} {x}}

\newcommand{\vv}{\vc{v}}

\newcommand{\vw}{\vc{w}}
\newcommand{\vn}{\vc{n}}

\newcommand{\bfU}{{\bf U}}
\newcommand{\bfV}{{\bf V}}

\newcommand{\FF}{\mathfrak F}

\newcommand{\vG}{\vc{G}}
\newcommand{\vphi}{{\boldsymbol \varphi}}
\newcommand{\tenS}{{\mathbb S}}

\newcommand{\pravp}{{\mathcal O}}
\newcommand{\wto}{\rightharpoonup}
\newcommand{\sym}{\operatorname{sym}}
\newcommand{\di}{{\rm d}}
\newcommand{\spam}{\operatorname{span}}
\newcommand{\esssup}{\operatorname{esssup}}

\makeatletter
\newsavebox{\@brx}
\newcommand{\llangle}[1][]{\savebox{\@brx}{\(\m@th{#1\langle}\)}%
  \mathopen{\copy\@brx\kern-0.5\wd\@brx\usebox{\@brx}}}
\newcommand{\rrangle}[1][]{\savebox{\@brx}{\(\m@th{#1\rangle}\)}%
  \mathclose{\copy\@brx\kern-0.5\wd\@brx\usebox{\@brx}}}
\makeatother

\begin{document}

\title*{Non-Newtonian compressible fluids with stochastic right hand side}
\author{Pavel Ludv\'ik \and V\'aclav M\'acha}
\institute{Pavel Ludv\'ik \at Faculaty of Science, Palack\'y University Olomouc, 17. listopadu 12, 771 46 Olomouc, Czech Republic, \email{pavel.ludvik@upol.cz}\\
V\'aclav M\'acha \at Institute of Mathematics of the Czech Academy of Sciences, \v Zitn\'a 25, 115 67 Praha 1, Czech Republic, \email{macha@math.cas.cz}}
\maketitle

\abstract
{Fluids with shear-rate-dependent viscosity form a special class of non-Newtonian fluids that play a crucial role in continuum dynamics. We consider a compressible barotropic flow of such fluids, governed by a generalized Navier-Stokes system. Due to the nonlinearity in the elliptic term, obtaining the existence of a weak solutions is challenging. To address this, we introduce the concept of a measure-valued solution, whose existence is established in a very general setting using an abstract theorem on the existence of a Young measure. We also discuss the proof of this key result.}

\section{Introduction}

The compressible non-Newtonian fluids are fluids whose viscosity is not constant and it may vary with the gradient of velocity and/or pressure. This allows us to capture certain physical phenomena, for example die swell, delayed die swell, shear thinning/thickening, and so on. The mathematical description of such fluids was deduced by M\'alek and Rajagopal \cite{MaRa}. 

The existence of a weak solution to such system without white noise is quite delicate as there is a lack of compactness of the velocity. The first papers in this direction are due to Mamontov \cite{mamontov1}, \cite{mamontov2} who considered an exponential growth of the viscosity together with an isothermal pressure ($p(\varrho) = a\varrho$ for some $a>0$). The existence of a weak solution was further treated by Feireisl, Liao, and M\'alek \cite{FeLiMa} who worked with a special form of the stress tensor, and further by Zhikhov and Pastukhova \cite{zhpa} for quite general setting -- unfortunately it turned out to be wrong. 

To overcome this problem, the notion of {\em dissipative weak solution} was introduced. Its existence was proven by Abbattiello, Feireisl and Novotn\'y \cite{AbFeNo} for isentropic pressure and by Basari\'c \cite{basaric} for an isothermal pressure. 

There are several reasons why to work with the stochastic force term in such system -- it can contain some numerical as well as physical uncertainties. This is the main aim of this article -- we introduce a notion of {\em measure-valued solution} for the nonlinear compressible Navier-Stokes equation with the stochastic force term and we prove the existence of both. 

Concerning the notion of the measure-valued solution for the stochastic equations we refer to the article by Hofmanov\'a, Koley and Sarkar \cite{HoKoSa} where this notion was given for the compressible Euler equations. Further, the stochastic version of the dissipative weak solution (martingale dissipative solution) was introduced by Hofmanov\'a, Zhu and Zhu in \cite{HoZhZh} who worked with the incompressible Euler system.

This paper is devoted to a nonlinear system of the form
\begin{equation}\label{main.sys}
\begin{split}
{\rm d}(\varrho \vu) + \diver(\varrho \vu \otimes \vu){\rm d}t + \diver \tenS {\rm d}t - \nabla p(\varrho){\rm d}t & = G(\varrho,\varrho\vu){\rm d}W\\
{\rm d} \varrho + \diver (\varrho \vu) {\rm d}t& = 0
\end{split}
\end{equation}
considered on $(0,T)\times\Omega$ where $T>0$ and $\Omega\subset \mathbb R^3$ is a bounded Lipschitz domain. The unknowns are the velocity $\vu:(0,T)\times \Omega \to \mathbb R^3$ and the density $\varrho:(0,T)\times \Omega\to [0,\infty)$. Furthermore, we assume that $\mathbb S\in \mathbb R^{3\times 3}_{sym}$ satisfies
\begin{equation}\label{young}
\tenS:D\vu = F^* (\tenS) + F(D\vu)
\end{equation}
for a suitable convex non-negative potential $F:\mathbb R^{3\times 3}_{sym}\to \mathbb R$ -- here $F^*$ denotes a convex conjugate of $F$ and $D\vu$ denotes the symmetric part of the gradient of $\vu$, i.e., $2D\vu = \left(\nabla \vu + (\nabla \vu)^T\right)$. We just recall that, assuming $F$ is a proper convex lower-semicontinuous \eqref{young} is equivalent to 
$$
\tenS\in \partial F(D\vu)
$$
or
$$
D\vu\in \partial F^*(\tenS)
$$
where $\partial$ denotes the sub-differential of $F$. Taking 
$$
F(D\vu) = \frac \mu2 |D\vu|^2 + \frac \lambda 2 |\diver \vu|^2
$$
one obtain a system describing the Newtonian fluid. Such system has been extensively studied -- the existence of a weak solution to a deterministic equation is provided by nowadays standard Lions-Feireisl theory (see \cite{lions} and \cite{FeNoPe}). Besides that, also the stochastic case has been examined -- we refer to book by Breit, Feireisl and Hofmanov\'a \cite{BrFeHo} and to a pioneering paper by Breit and Hofmanov\'a \cite{BrHo}. There are way more results dealing with the analysis of the Newtonian fluid and it is impossible to present an exhaustive list thus we just presented some of the most significant representatives.

We would also like to point out that the above implicit relation between $\tenS$ and $D\vu$ contains a lot of models, for example the power law fluids, i.e. fluids whose stress tensor satisfies for example
$$
\tenS(D\vu) = \left(1 + |D\vu|^2\right)^{(p-2)/2} D\vu
$$
which is attained by the choice $$F(D\vu) = \frac 1p \left(1+|D\vu|^2\right)^{p/2}.$$

The system \eqref{main.sys} is complemented with the homogenous Dirichlet boundary condition
\begin{equation}\label{dirichlet}
\vu{\restriction_{\partial\Omega}} = 0,\ \Omega\mbox{ is of class }C^{2+\nu},\, \nu>0
\end{equation}
which captures the case when the fluid sticks to the boundary of $\Omega$. 
\medskip \\
\indent
Let us specify the assumptions on data. We start with the definition of the stochastic forcing term. Let $(\mathcal O, \FF, (\FF_t)_{t\geq0},\mathbb P)$ be a stochastic basis with a complete right-continuous filtration and
$W$ be a cylindrical $(\FF_t)$-Wiener process in a separable Hilbert space $\mathfrak U$ formally given by
$$
W = \sum_{k=1}^\infty e_k W_k(t),
$$
where $\{W_k\}_{k=1}^{\infty}$ is a sequence of mutually independent real-valued Brownian motions relative to $(\FF_t)_{t\geq0}$ and $\{e_k\}_{k=1}^{\infty}$ is an orthonormal basis of $\mathfrak U$.

Next, we assume that 
\begin{equation}\label{G.ass.1}
G(\varrho,\vu) \, \di W = \sum_{k=1}^\infty {\bf G}_k(\varrho,\vu)\, \di W_k,\quad {\bf G}_k(\varrho,  \vu) = \varrho {\bf F}_k(\varrho,\vu)
\end{equation}
where 
\begin{equation*}
{\bf F}_k:\Omega\times [0,\infty)\times \mathbb R^3\to \mathbb R^3,  {\bf F}_k\in C^1(\Omega\times [0,\infty)\times \mathbb R^3)\end{equation*}
satisfies
\begin{equation*}
\|{\bf F}_k\|_{L^\infty(\Omega\times[0,\infty)\times \mathbb R^3)} + \|\nabla_\vu {\bf F}_k\|_{L^\infty(\Omega\times [0,\infty)\times\mathbb R^3)}\leq f_k,\quad \sum_{k=1}^\infty f_k^2<\infty
\end{equation*}
and there is some positive constant $\alpha$ such that ${\bf F}_k(\varrho,\vu)$ is zero whenever $\varrho<\alpha$, $\varrho >\frac 1\alpha$, or $|\vu|>\frac 1\alpha$ -- this last assumption is here just for simplicity as there is a known method allowing to remove it (we refer to \cite{BrFeHo}).  We define $\mathbb F = \{{\bf F}_k\}_{k=1}^\infty$.

The initial condition is $\varrho(0,\cdot) = \varrho_0(\cdot)$ and $\varrho(0,\cdot)\vu(0.\cdot) = (\varrho \vu)_0(\cdot)$ and the initial data are supposed to satisfy
\begin{equation*}
\int_\Omega P(\varrho_0)\, \di x <\infty,\quad \int_\Omega\frac 12 \frac{(\varrho \vu)_0^2}{\varrho_0}\, \di x <\infty,\ \mbox{and }\int_\Omega (\varrho\vu)_0\, \di x <\infty
\end{equation*}
$\mathbb P$-almost surely.
Here $P$ stands for the pressure potential introduced below \eqref{yes.phi.energy}. However, this initial datum is given only by its law -- we refer to \eqref{lambda.ini.con.1} and \eqref{lambda.ini.con.2} for the precise statement of the initial condition.

Next, the potential $F$ is supposed to be of class $C^1$, to satisfy $F(0)=0$, and we assume the existence of an $N$-function $g:[0,\infty)\to [0,\infty)$ (see \cite[Chapter 8]{AdFo}) such that
\begin{equation}\label{F.potential.2}
F(D\vu) \geq g(|D\vu|).
\end{equation}
\bigskip \\
\indent Up to our knowledge the generalized compressible Navier-Stokes equation with a stochastic right hand side as presented above has not been studied yet.

The main aim of this paper is to provide a notion of a solution which is on one hand sufficiently general so one can still deduce its existence and, on the other hand, it still posses the weak-strong uniqueness property. 

To achieve this, we use a theorem on the existence of Young measures for a sequence of function-valued random variables, as introduced in \cite{BrFeHo}. We revisit this theorem under relaxed assumptions and provide a detailed proof in Section \ref{m-v.s.a.i.e}.

Section \ref{S.t.G.A} focuses on the analysis of the Galerkin approximation, which serves as the starting point for the subsequent limit process, leading to the existence of a properly defined measure-valued solution to \eqref{main.sys}. The existence of solution to this approximation is established using the methods presented in \cite{BrFeHo}.

The limit in the Galerkin approximation is presented in Section \ref{m-v.m.s}, which also provides the proper definition of the measure-valued solution. The existence of such a solution is summarized in this section as well, stated in a theorem that we consider the main result of this paper.


\section{Remarks on Young measures} \label{m-v.s.a.i.e}

In our analysis below, we make use of two standard theorems that appear in the literature in various forms. We state the specific variants needed for our purposes here.

\begin{Theorem}[{de la Vallée Poussin criterion, see \cite[Theorem 4.5.9.]{Bo}}]\label{t:poussin}
Let $(X,\mathcal A,\mu)$ be a measure space with a finite measure. A family $\mathcal F$ of $\mu$-integrable functions is equi-integrable if and only if there exists a non-decreasing function $g:[0,+\infty)\to[0,+\infty)$ such that
\[
\lim_{t\to+\infty}\frac{g(t)}{t}=\infty \mbox{ and } \sup_{f\in\mathcal F} \int_{X} g(|f(x)|)\,\rm{d}\mu(x)\leq\infty.
\]
In such a case, one can choose a convex non-decreasing function $g$.
\end{Theorem}

\begin{Theorem}[{Dunford-Pettis theorem, see \cite[Theorem 4.7.18.]{Bo}}]\label{t:du-pe}
Let $(X,\mathcal A,\mu)$ be a measure space with a finite measure and let $\mathcal F$ be a family of $\mu$-measurable functions. The set $\mathcal F$ is equi-integrable precisely when it is weakly relatively compact in $L^1(X)$.
\end{Theorem}

The main tool yielding the measure-valued solution is the concept of the Young measures (see c.f. \cite[Theorem 6.2]{Pedregal}) adapted to the stochastic setting. The theorem following can be found as Lemma 2.2 in \cite{HoKoSa}:

\begin{Lemma}
Let $N,M\in \mathbb N,\, \Omega\subset \mathbb R^N$, let $(\pravp,\mathfrak F,\mathbb P)$ be a complete probability space and $\{\bfU_n\}_{n=1}^{\infty},\, \bfU_n:\Omega\times\pravp\to \mathbb R^M$, be a sequence of random variables such that\footnote{Hereinafter $C$ is reserved for a constant which may vary from line to line but remains independent of solutions.  Any specific dependencies will be explicitly highlighted when necessary.}
$$
\mathbb E\left[\|\bfU_n\|_{L^p(\Omega)}^p\right]\leq C,\, \mbox{for a certain }p\in (1,\infty).
$$
Then there exists a new subsequence $\{\tilde{\bfU}_n\}_{n=1}^{\infty}$ (not relabeled), defined on the standard probability space $\left([0,1],\overline{\mathcal B[0,1]},\mathbb P\right)$, and a parametrized family $\left\{\tilde{\nu}_{\omega, x}\right\}_{(\omega,x)\in \pravp\times\Omega}$ of random probability measures on $\mathbb R^M$, regarded as a random variable taking values in $(L^\infty_{\omega^*}(\Omega;\mathcal P(\mathbb R^M)),w^*)$, such that $\bfU_n$ has the same law as $\tilde \bfU_n$, i.e., $\bfU_n \sim_d\tilde\bfU_n$, and the following property holds: for any Carath\'eodory function $G = G(x,Z)$, $x\in \Omega$, $Z\in \mathbb R^M$, such that
$$
|G(x,Z)|\leq C(1 + |Z|^q),\ 1\leq q<p,\ \mbox{uniformly in }x,
$$
implies $\mathcal L$-a.s.
$$
G(\cdot,\tilde \bfU_n)\rightharpoonup \tilde G\ \mbox{in }L^{p/q}(\Omega)
$$
where
$$
\overline G(x) = \langle \tilde{\nu}_{\omega,x}; G(x,\cdot)\rangle:=\int_{\mathbb R^M} G(x,z)\, {\rm d}\tilde{\nu}_{\omega,x}(z)
$$
for almost all $x\in \Omega$.
\end{Lemma}

Here, we present a stronger theorem together with its proof. Namely
\begin{Theorem}\label{thm:hofmanova1}
Let $N,M\in \mathbb N$, $\Omega\subset\mathbb R^N$ of finite Lebesgue measure, let $(\pravp,\mathfrak F,\mathbb P)$ be a complete probability space and $\{\bfU_n\}_{n=1}^{\infty},\, \bfU_n:\Omega\times\pravp\to \mathbb R^M$, be a sequence of random variables such that
\begin{equation*}
\sup_{n\in\mathbb N} \int_{\mathcal O}\int_{\Omega} g(|\bfU_n|) {\rm d}\lambda(x) {\rm d}\mathbb{P}(\omega)\leq C,
\end{equation*}
for some $g:[0,\infty)\to[0,\infty)$ which is a non-decreasing function such that $\lim_{t\to\infty}\frac{g(t)}{t}=\infty$. 

Then there exists a subsequence $\{\tilde{\bfU}_n\}_{n=1}^\infty$ and a parametrized family $\left\{\tilde{\nu}_{\omega,x}\right\}_{(\omega,x)\in \pravp\times\Omega}$ as in the previous lemma such that for every Carath\'eodory function $G = G(x,Z)$ such
$$
G(x,\tilde{\bfU}_n) \ \mbox{is weakly convergent in}\ L^1(\Omega)\ \mbox{almost surely}
$$
it holds that
$$
G(x,\tilde \bfU_n) \rightharpoonup \langle \tilde{\nu}_{\omega,x}; G(x,\cdot)\rangle
$$
in $L^1(\Omega)$ almost surely.
\end{Theorem}

\begin{proof}
Let $(\nu_n)_{\omega,x}:= \delta_{\bfU(\omega,x)}$ and let $\mathcal L[\nu_n]$ be its law. Next, let 
$$
B_{R,1} = \left\{\nu\in L^\infty(\Omega;\mathcal P(\mathbb R^M),w^*),\, \int_X\int_\mathbb R |\xi| \, {\rm d}\nu_x(\xi){\rm d}\lambda(x)\leq R\right\}.
$$
This set is relatively compact in $L^\infty(\Omega;\mathcal P(\mathbb R^M), w^*)$. This can be justified by the observation that in \cite[Theorem 5]{DeVo}, without changing the proof, we can replace $\mathbb R$ by $\mathbb R^M$ and the assumption that $(X,\lambda)$ is a finite measure space with the assumptions that $(X,\lambda)$ is a $\sigma$-finite space and the corresponding $L^1(X)$ is separable. Consequently, by the Chebyshev inequality
\begin{multline*}
\mathcal L[\nu_n](B_{R,1}^c) = \mathbb P\left(\nu_n\notin B_{R,1}\right) = \mathbb P\left(\int_{\Omega}|\bfU_n|\, {\rm d}x>R\right)\\
\leq \frac 1R \mathbb E\left[\|\bfU_n\|_{L^1(\Omega)}\right] \leq \frac CR.
\end{multline*}
The family of laws $\{\mathcal L[\nu_n]\}_{n=1}^{\infty}$ is therefore tight in 
$L^\infty(\Omega;\mathcal P(\mathbb R^M), w^*)$.

Next, we show that also the family of laws $\{\mathcal L[\bfU_n]\}_{n=1}^\infty$ is tight on $L^1(\Omega,w)$. For the fixed $R>0$, we consider the set 
\[
\left\{\bfV\in L^1(\Omega),\, \|g(\bfV)\|_{L^1(\Omega)}\leq R\right\},
\] 
which is equi-integrable in $L^1(\Omega,w)$ by Theorem~\ref{t:poussin}. Thus, it is also relatively weakly compact by Theorem~\ref{t:du-pe}. 

Therefore, 
$$
\mathcal L[\bfU_n]\left(\left\{\bfV\in L^1(\Omega),\,\|g(\bfV)\|_{L^1(\Omega)}>R\right\}\right)\leq  \frac 1{R} \mathbb E\left[\|g(\bfU_n)\|_{L^1(\Omega)}\right]\leq \frac C{R},
$$
by Chebyshev inequality, and we thus obtain the desired tightness.

Thus, the family of joint laws $\{\mathcal L[\bfU_n,\nu_n],\, n\in\mathbb{N}\}$ is tight in $L^1(\Omega,w)\times L^\infty(\Omega;\mathcal P(\mathbb R^M), w^*)$ and we may apply Jakubowski-Skorokhod's theorem (see \cite[Theorem 2]{Ja} or \cite[Theorem 2.7.1]{BrFeHo}) to obtain a complete probability space $(\tilde{\pravp},\tilde{\mathcal{S}},\tilde{\mathbb{P}})$ and Borel random variables $[\tilde{\bfU}_n,\tilde{\nu}_n]$ with the law $\mathcal L[\bfU_n,\nu_n]$, for $n\in\mathbb{N}$, and converging $\tilde{\mathbb P}$-a.s. to a Borel random variable $[\tilde{\bfU},\tilde{\nu}]$ possessing values in $L^1(\Omega,w)\times L^\infty(\Omega;\mathcal P(\mathbb R^M), w^*)$.

Now, let us define
\begin{multline*}
M := \left\{(\bfV,\nu)\in L^1(\Omega)\times L^\infty(\Omega;\mathcal P(\mathbb R^M), w^*),\, \int_{\Omega} \psi(x) \int_{\mathbb R^M} \phi(\xi){\rm d}\nu_x(\xi){\rm d}x \right.\\
\left.= \int_{\Omega} \psi(x)\phi(\bfV){\rm d}x \mbox{ for all $\psi\in L^1(\Omega)$, $\phi\in C_b(\mathbb R^M)$}\right\}.
\end{multline*}
Then, from the definition of $\nu_n$, $n\in\mathbb{N}$, and the equality of laws of $[\bfU_n,\nu_n]$ and $[\tilde{\bfU}_n,\tilde{\nu}_n]$ follows
\begin{multline*}
1 = \mathbb P\left(\int_{\Omega} \psi(x) \int_{\mathbb R} \phi(\xi){\rm d}\nu_{n,x}(\xi){\rm d}x =\right.\\ \left.
\int_{\Omega} \psi(x)\phi(\bfU_n){\rm d}x \mbox{ for all $\psi\in L^1(\Omega)$, $\phi\in \mathcal C_b(\mathbb R^M)$}\right) \\
= \mathcal L[\bfU_n,\nu_n](M) = \mathcal L[\tilde{\bfU}_n,\tilde{\nu}_n](M) \\
=\mathbb {\tilde P}\left(\int_{\Omega} \psi(x) \int_{\mathbb R^M} \phi(\xi){\rm d}\tilde{\nu}_{n,x}(\xi){\rm d}x \right.\\ \left.
=\int_{\Omega} \psi(x)\phi(\tilde\bfU_n){\rm d}x \mbox{ for all $\psi\in L^1(\Omega)$, $\phi\in \mathcal C_b(\mathbb R^M)$}\right).
\end{multline*}
We can conclude that $\tilde \nu_{n,\omega,x}(\cdot)=\delta_{\hat\bfU_n(\omega,x)}(\cdot)$ for all $\omega$ from a set of full probability measure $\tilde {\mathbb P}$.

Finally, let us take $\omega$ from a subset of $\pravp$ of full measure such that $\tilde \nu_{n,\omega,x}(\cdot)=\delta_{\hat\bfU_n(\omega,x)}(\cdot)$ together with the weak convergence of $[\tilde\bfU_n,\tilde\nu_n]$ to $[\tilde\bfU,\tilde\nu]$. If we denote 
\[
\overline{\phi(\tilde\bfU)}:=\langle \tilde{\nu}_{\omega,x}; \phi(\cdot)\rangle = \int_{\mathbb R^M} \phi(\xi){\rm d}\tilde{\nu}_{\omega,x}(\xi){\rm d}x,
\]
we have, for all $\psi\in L^1(\Omega)$, $\phi\in \mathcal C_b(\mathbb R^M)$,
\begin{align*}
\lim_{n\to\infty} \int_{\Omega} \psi(x)\phi(\tilde\bfU_n){\rm d}x &= \lim_{n\to\infty} \int_{\Omega} \psi(x) \int_{\mathbb R^M} \phi(\xi){\rm d}\tilde{\nu}_{n,x}(\xi){\rm d}x \\
&= \int_{\Omega} \psi(x) \int_{\mathbb R^M} \phi(\xi){\rm d}\tilde{\nu}_{x}(\xi){\rm d}x
= \int_{\Omega} \psi(x) \overline{\phi(\tilde\bfU)}{\rm d}x.
\end{align*}

We can conclude that $\phi(\tilde\bfU_n)\rightharpoonup\overline{\phi(\tilde\bfU)}$ in $L^1(\Omega,w)$ for every $\phi\in\mathcal C_b(\mathbb R^M)$.

  On the other hand, we can apply the theorem stating the existence of a Young measure (see \cite[Theorem 6.2]{Pedregal} or \cite[Theorem 2.8.2.]{BrFeHo}). Since the probability measures are determined by their action on continuous bounded functions, the obtained Young measure must coincide with $\tilde{\nu}_{\omega,x}$ for every $\omega\in\Omega$. Knowing this, we can extend our previous conclusion from $\phi\in\mathcal C_b(\mathbb R^M)$ to any 
Carath\'eodory function $G = G(x,Z)$ such $G(x,\tilde{\bfU}_n)$ is weakly convergent in $L^1(\Omega)$ almost surely. This gives the demanded claim.
\end{proof}

As stated in \cite[Lemma 2.3]{HoKoSa}, by relaxing our requirements, we can still obtain a reasonable outcome in terms of the concentration defect measures. We present the lemma along with its proof (which is omitted in the cited reference), which is a stochastic modification of the proof of \cite[Lemma 2.1]{FeGwSwWi}.

Let us recall that a sequence of functions generating a Young measure (i.e. satisfying the conclusion of Theorem \ref{thm:hofmanova1}) is weakly convergent in $L^1$ (see e.g. \cite[Remark 3.131]{FlGo}).

\begin{Lemma}
Let $N,M\in \mathbb N$, $\Omega\subset\mathbb R^N$ of finite Lebesgue measure, let $(\pravp,\mathfrak F,\mathbb P)$ be a complete probability space and $\{\bfU_n\}_{n=1}^{\infty},\, \bfU_n:\Omega\times\pravp\to \mathbb R^M$, be a sequence of random variables generating a Young measure $\{\nu_{\omega,x}\}_{x\in\Omega}$. Let $G:\mathbb R^M\to[0,\infty)$ be a continuous function such that 
\[
\sup_{n\in\mathbb N}\mathbb E \left[\|G(\bfU_n)\|^p_{L^1(\Omega)}\right]<\infty, \mbox{ for a certain $p\in(1,\infty)$},
\]
and let $F$ be continuous function satisfying
\[
F:\mathbb R^M\to \mathbb R, \quad |F(z)| \leq G(z), \mbox{ for all $z\in\mathbb R^M$.}
\]
Let us denote $\mathbb P$-a.s.
\[
F_\infty := \tilde F - \langle \nu_{\omega,x};F(\cdot)\rangle\dx, \quad G_\infty := \tilde G - \langle \nu_{\omega,x};G(\cdot)\rangle\dx
\]
where $\tilde{F},\tilde{G}\in\mathcal M(\Omega)$ (a space of bounded Radon measures on $\Omega$) are $w^*$-limits of $\{F(\bf U_n)\}_{n=1}^\infty$, $\{G(\bf U_n)\}_{n=1}^\infty$ respectively in $\mathcal M(\Omega)$. Then $|F_\infty|\leq G_\infty$ $\mathbb P$-a.s.
\end{Lemma}

\begin{proof}
First, we express $\tilde{F}$ and $\tilde{G}$ as the $w^*$-limits of the sequences $\{F(\mathbf{U}_n)\}_{n=1}^\infty$ and $\{G(\mathbf{U}_n)\}_{n=1}^\infty$ in $\mathcal M(\Omega)$, respectively, while also incorporating the truncation parameter $M > 0$. Hence,
$$
\begin{aligned}
\langle\tilde F; \phi\rangle&=\lim _{n \to \infty} \int_{\Omega \cap\{\left|\bfU_n\right| \leq M\}} F\left(\bfU_n\right) \phi \dx +\lim _{n \to \infty} \int_{\Omega \cap\{\left|\bfU_n\right| > M\}} F\left(\bfU_{n}\right) \phi \dx, \\
 \langle\tilde G; \phi\rangle &= \lim_{n \to \infty} \int_{\Omega \cap\{\left|\bfU_n\right| \leq M\}} G\left(\bfU_n\right) \phi \dx +\lim_{n \to \infty} \int_{\Omega \cap\{\left|\bfU_n\right| > M\}} G\left(\bfU_n\right) \phi \dx,
\end{aligned}
$$
for every $\phi\in\mathcal C_b(\Omega)$ $\mathbb P$-a.s.

Because $F,G$ are bounded continuous on the set $\{z\in\Omega,\, |z|\leq M\}$, we can employ a Young measure $\{\nu_{\omega,x}\}_{x\in\Omega}$ to obtain
$$
\begin{aligned}
&\lim _{n \to \infty} \int_{\Omega \cap\{\left|\bfU_n\right| \leq M\}} F\left(\bfU_{n}\right) \phi \dx = \int_{\Omega \cap\{\left|\bfU_n\right| \leq M\}}\left\langle \nu_{\omega,x};F(\cdot)\right\rangle\phi\dx, \\
&\lim_{n \to \infty} \int_{\Omega \cap\{\left|\bfU_n\right| \leq M\}} G\left(\bfU_n\right) \phi \dx = \int_{\Omega \cap\{\left|\bfU_n\right| \leq M\}} \left\langle \nu_{\omega,x};G(\cdot)\right\rangle\phi\dx,
\end{aligned}
$$
for every $\phi\in\mathcal C_b(\Omega)$ $\mathbb P$-a.s.

Then we get, for almost all $\omega\in\mathcal O$ and every $\phi\in\mathcal C_b(\Omega)$,
$$
\begin{aligned}
\lim_{M \to \infty}\left(\lim _{n \to \infty} \int_{\Omega \cap\{\left|\bfU_n\right| \leq M\}} F\left(\bfU_n\right) \phi \dx\right)=\int_{\Omega}\left\langle \nu_{\omega,x} ; F(\cdot)\right\rangle \phi\dx, \\
\lim_{M \to \infty}\left(\lim_{n \rightarrow \infty} \int_{\Omega \cap\{\left|\bfU_n\right| \leq M\}} G\left(\bfU_n\right) \phi \dx\right)=\int_{\Omega}\left\langle\nu_{\omega,x} ; G(\cdot)\right\rangle \phi\dx .
\end{aligned}
$$

By combining the above considerations, we arrive at the conclusion that
$$
\begin{aligned}
\left\langle F_{\infty} ; \phi\right\rangle & =\lim _{M \to \infty}\left(\lim _{n \to \infty} \int_{\Omega \cap\{\left|\bfU_n\right| > M\}} F\left(\bfU_n\right) \phi \dx\right),\\
\left\langle G_{\infty};\phi\right\rangle & =\lim _{M \to \infty}\left(\lim _{n \to \infty} \int_{\Omega \cap\{\left|\bfU_n\right| > M\}} G\left(\bfU_n\right) \phi \dx\right) ,
\end{aligned}
$$
for every $\phi\in\mathcal C_b(\Omega)$ $\mathbb P$-a.s.

Since $|F| \leq G$, it follows that 
$$
\left|\left\langle F_{\infty} ; \phi\right\rangle\right| \leq \left\langle G_{\infty};\phi\right\rangle,
$$
for every $\phi\in\mathcal C_b(\Omega)$ $\mathbb P$-a.s. which completes the proof of $|F_\infty|\leq G_\infty$ $\mathbb P$-a.s.
\end{proof}

In what follows, we give the proof of the existence of a measure-valued solution to \eqref{main.sys} endowed with the Dirichlet boundary condition \eqref{dirichlet}. We do not expect any complication when treating other types of well known boundary conditions like the full-slip boundary condition or the periodic boundary condition.

\section{Galerkin approximation}\label{S.t.G.A}
Since $L^2(\Omega)$ has a countable Schauder basis consisting of smooth functions with zero on the boundary, we consider $n$-dimensional spaces $X_n\subset C^\infty_0(\Omega)$ with the property 
\begin{equation}\label{hustota.v.L2}L^2(\Omega) = \overline{\bigcup_n X_n}^{\|\cdot\|_{L^2}}.\end{equation}
Let $\Pi_n$ be the orthonormal projection $\Pi_n:L^2(\Omega)\to X_n$. In this section we tackle the problem of the existence of a pair $(\varrho_n,\vu_n)\subset C([0,T], L^1)\times C([0,T], X_n)$ satisfying the system
\begin{equation}
\begin{split}\label{approximace}
{\rm d}\varrho_n + \diver (\varrho_n \vu_n)\, {\rm d}t & = 0 \\
{\rm d}\Pi_n(\varrho_n\vu_n) + \Pi_n\big((\diver \varrho_n\vu_n\otimes \vu_n) +  \nabla p(\varrho_n)\big)\, {\rm d}t& =  \Pi_n\big(\diver \mathbb S_n\big)\, {\rm d}t\\
\ + \Pi_n(&\varrho\Pi_n(\mathbb F(\varrho,\vu)))\, {\rm d}W.
\end{split}
\end{equation}

\begin{Definition}\label{Galerkin.definition}
Let $\Lambda$ be a Borel probability measure on $C^{2+\nu}(\Omega)\times X^n$. Then $((\mathcal O, \mathfrak F, (\mathfrak F_t)_{t\geq 0},\mathbb P),\varrho_n, \vu_n,W)$ is called a martingale solution to \eqref{approximace} with initial datum $\Lambda$ if the following is satisfied:
\begin{enumerate}
\item $(\mathcal O, \mathfrak F, (\mathfrak F_t)_{t\geq 0},\mathbb P)$ is a stochastic basis with a complete right-continuous filtration;
\item $W$ is a cylindrical Wiener process;
\item the density $\varrho_n$ is $(\mathfrak F_t)$-progressively measurable and satisfies
$$
\varrho \in L^p((0,T)\times \Omega)\, \mbox{ for any }p\in (1,\infty) ,\ \varrho_n\geq \underline \varrho >0\quad \mathbb P\mbox{-a.s.};
$$
\item the velocity $\vu_n$ is $(\mathfrak F_t)$-progressively measurable and satisfies
$$
\vu_n \in C([0,T],\, X_n)\quad \mathbb P\mbox{-a.s.};
$$
\item there exists an $\mathfrak F_0$ measurable random variable $(\varrho_0,\vu_0)$ such that  
$$
\mathcal L(\varrho_n(0), \vu_n(0)) = \Lambda;
$$
\item the approximate continuity equation
\begin{equation*}
\pat \varrho_n + \diver(\varrho_n\vu_n) = 0
\end{equation*}
holds in $(0,T)\times \Omega$ $\mathbb P$-a.s. and $\varrho(0) = \varrho_0$ $\mathbb P$-a.s.;
\item the approximate momentum equation
\begin{multline*}
-\int_0^T\partial_t\Phi\int_\Omega \varrho_n\vu_n\varphi \, {\rm d}x{\rm d}t - \Phi(0)\int_\Omega \varrho_0\vu_0 \cdot \varphi \ {\rm d}x \\
= \int_0^T\Phi \int_{\Omega} \left(\varrho_n \vu_n \otimes \vu_n:\nabla\varphi + p(\varrho_n)\diver \varphi\right)\, {\rm d}x{\rm d}t\\
 - \int_0^T\Phi\int_\Omega \mathbb S(\nabla \vu_n):\nabla\varphi \, {\rm d}x{\rm d}t\\ + \int_0^T\Phi \int_\Omega \varrho_n \Pi_n[\mathbb F(\varrho_n,\vu_n)]\cdot \varphi \, {\rm d}x{\rm d}W
\end{multline*}
holds for all $\Phi\in C^\infty_c([0,T))$ and all $\varphi \in X_n$ $\mathbb P$-a.s.;
\item the energy inequality
\begin{multline}\label{yes.phi.energy}
-\int_0^T\pat \Phi \int_\Omega \left(\frac 12 \varrho |\vu|^2 + P(\varrho)\right)\, \di x \di t - \Phi(0)\int_\Omega \left(\frac 12 \varrho |\vu|^2 + P(\varrho)\right)\, \di x\\
 + \int_0^T\Phi \int_\Omega \mathbb S(D\vu):\nabla\vu\, \di x \di t = \frac 12 \sum_{k=1}^\infty \int_0^T\Phi\int_\Omega \varrho \left|\Pi_n ({\bf F}_k(\varrho,\vu))\right|^2\, \di x \di t\\
 \int_0^T \Phi \int_\Omega \varrho \Pi_n \left[\mathbb F(\varrho,\vu)\right]\cdot \vu \, \di x \di W
\end{multline}
 holds for all $\Phi \in C^\infty_c([0,T))$ $\mathbb P$-a.s.
 Here $P$ stands for the usual pressure potential, i.e. 
$$
P(\varrho) = \varrho\int_1^\varrho\frac{p(z)}{z^2}\, {\rm d}z
$$
\end{enumerate}
and $\mathbb S = \mathbb S(D\vu_n)$.
\end{Definition}
Let us mention that the energy inequality and the Burkholder-Davis-Gundy inequality yield
\begin{multline}\label{no.phi.energy}
\mathbb E \left[\left|\sup_{\tau\in [0,T]} \int_\Omega \frac 12 \varrho |\vu|^2 + P(\varrho)\, \di x \right|^r\right] 
+ \mathbb E\left[\left|\int_0^T\int_\Omega \mathbb S:\nabla\vu\, \di x\di t\right|^r\right]\\
\leq C \mathbb E\left[\left|\int_\Omega \frac 12 \varrho_0 |\vu_0|^2 + P(\varrho_0)\, \di x \right|^r + 1\right]
\end{multline}
for any $r>1$.

The starting point of the existence of the measure-valued solution is the following theorem. Its proof is recalled for reader's convenience in the next section.
\begin{Theorem}\label{Galerkin.existence}
Let $\Lambda$ be a Borel probability measure on $C^{2+\nu}(\Omega)\times X^n$ such that 
$$
\Lambda [0<\underline \varrho \leq \varrho, \| \varrho\|_{C^{2+\nu}(\Omega)}\leq \overline \varrho] = 1,\quad \int \|\vv\|_{X^n}^r\, \di \Lambda [\varrho,\vv]\leq \overline u
$$
for some deterministic constants $\underline \varrho,\, \overline \varrho,\, \overline u$ and some $r>2$. Then there exists the martingale solution specified in Definition \ref{Galerkin.definition}.
\end{Theorem}

\subsection{Proof of Theorem \ref{Galerkin.existence}}

Let $\{{\mathbf \omega_k}\}_{k=1}^\infty$ be the orthogonal basis of $W^{3,2}_0(\Omega)$ formed by the smooth eigenvectors of the problem
\begin{equation*}
\begin{split}
-\Delta^3 \mathbf \omega = \lambda {\mathbf \omega},\\
{\mathbf \omega}\restriction_{\partial\Omega} = 0.
\end{split}
\end{equation*}
(see \cite[Theorem 4.11 \& Remark 4.14 in the Appendix]{MNRR} for the existence of such basis). We define $X_n = \spam\{\omega_k, k=1,\ldots,n\}$ as the sequence of spaces satisfying \eqref{hustota.v.L2}.
\bigskip \\
\indent Fix $n\in \mathbb N$. In the first step, we show the solvability of the system consisting of 
\begin{itemize}
\item the continuity equation
\begin{equation}\label{proof.compressible1}
\pat \varrho  + \diver (\varrho [\vu]_R) = \varepsilon \Delta \varrho
\end{equation}
endowed with the boundary condition 
$$
\frac{\partial \varrho}{\partial \vn} = 0\, \mbox{ on }\partial \Omega;
$$
\item and the momentum equation
\begin{multline}\label{proof.momentum1}
\Pi_n[\di (\varrho \vu)] + \Pi_n[\diver (\varrho [\vu]_R \otimes\vu)]\,\di t - \mu \Pi_n[\Delta^3 \vu]\, \di t\\ - \Pi_n[\diver \mathbb S(D\vu)]\, \di t
+ \Pi_n[\chi (\|\vu_n - R\|_{X_n})\nabla p(\varrho)]\, \di t\\ - \varepsilon\Pi_n[\Delta(\varrho \vu)]\, \di t = \Pi_n[\varrho \Pi_n[\mathbb F(\varrho,\varrho \vu)]]\,\di W,
\end{multline}
\end{itemize}
where $\mu>0$, $\varepsilon>0$, and $R>0$ are given parameter.  The meaning of $[\cdot]_R$ and $\chi$ can be found in \cite[Chapter 4]{BrFeHo}. We look for the solution in the form 
$$ \vu = \sum_{i=1}^N c_i {\mathbf \omega}_i,$$
where $c_i\in \mathbb R$ are unknowns assuming the stochastic basis $(\Omega,\mathfrak F, (\mathfrak F_t)_{t\geq 0}, \mathbb P)$ and the cylindrical process $W$ are given. The system \eqref{proof.compressible1} and \eqref{proof.momentum1} is uniquely solvable for initial conditions satisfying  
\begin{multline*}
\varrho_0\in C^{2+\nu}(\Omega),\quad \vu_0\in X_n,\quad \mathbb P\{0<\underline \varrho\leq \varrho_0,\, \|\varrho_0\|_{C^{2+\nu}}\leq \overline \varrho \} = 1,\\ \quad \mathbb E[\|\vu_0\|_{X_n}^r]\leq \overline u,
\end{multline*}
where $\underline \varrho,\, \overline \varrho,\, \overline u,\, r >2$. Moreover, we choose $\varrho_0$ and $\vu_0$ such that $\mathcal L(\varrho_0,\vu_0) = \Lambda$ and $\varrho_0$ and $\vu_0$ are $\mathfrak F_0$ measurable. 
The above mentioned existence of the solution is shown (with slight modification) in \cite[Theorem 4.1.12]{BrFeHo}.
In the rest of the proof we pass with parameters $\mu$ and $\varepsilon$ to zero and with $R$ to infinity. \medskip \\
\indent
First we tend with $\mu$ to zero. The energy balance
\begin{multline*}
-\int_0^T\pat \Phi\int_\Omega \frac 12 \varrho |\vu|^2 + P(\varrho)\, \di x \di t - \Phi(0) \int_\Omega \frac 12 \varrho |\vu|^2 + P(\varrho)\, \di x\\
 + \int_0^T\Phi\int_\Omega \mathbb S(D\vu):\nabla \vu + \mu|\nabla^3\vu|^2 + \varepsilon \varrho |\nabla \vu|^2 + \varepsilon P''(\varrho)|\nabla \varrho|^2\, \di x \di t\\
   = \frac 12 \sum_{k=1}^\infty \int_0^T\Phi \int_\Omega \varrho |\Pi_n[F_k(\varrho,\vu)]|^2\, \di x \di t\\ + \int_0^T\Phi\int_\Omega \varrho \Pi_n [\mathbb F(\varrho,\vu)]\cdot \vu \, \di t \di W
\end{multline*}
holds for all $\Phi\in C_c^\infty([0,T))$ $\mathbb P$-almost surely -- see \cite[Proposition 4.1.14]{BrFeHo}. The Burkholder-Davis-Gundy inequality and the Gronwall lemma yield
\begin{multline}\label{proof.energy2}
\mathbb E \left[\left|\sup_{\tau\in [0,T]} \int_\Omega \frac 12 \varrho |\vu|^2 + P(\varrho)\, \di x \right|^r\right] \\
+ \mathbb E\left[\left|\int_0^T\int_\Omega \mathbb S:\nabla\vu + \varepsilon \varrho |\nabla\vu|^2 + \mu |\nabla^3\vu|^2 + \varepsilon P'' (\varrho)|\nabla\varrho|^2\, \di x\di t\right|^r\right]\\
\leq C\, \mathbb E\left[\left|\int_\Omega \frac 12 \varrho_0 |\vu_0|^2 + P(\varrho_0)\, \di x \right|^r + 1\right],
\end{multline}
whenever $r\geq 2$. This (together with \eqref{young}) enables us to deduce the following estimates
\begin{equation}
\begin{split}\label{proof.estimates.1}
\mathbb E \left[\left|\int_0^T\int_\Omega F(D\vu)\, \di x\di t\right|^r\right]&\leq C,\\
\mathbb E \left[ \left\|\varrho |\vu|^2\right\|_{L^\infty((0,T);L^1(\Omega))}^r\right]&\leq C,\\
\mathbb E\left[\left\| P (\varrho)\right\|_{L^\infty((0,T);L^1(\Omega))}^r\right]&\leq C.
\end{split}
\end{equation}
Since $[\vu]_R$ is bounded, \eqref{proof.compressible1} yields the existence of $C>0$ such that 
\begin{equation}\label{proof.odhad.rho1}
\begin{split}
\esssup_{t\in [0,T]} \left(\|\varrho(t)\|_{C^{2+\nu}(\Omega)} + \|\pat \varrho\|_{C^\nu(\Omega)} \right)\leq C,  \\
\frac 1C \leq \varrho(x)\leq C,\quad \mbox{for all }(t,x)\in [0,T)\times \Omega.
\end{split}
\end{equation}
This can be found for example at \cite{NoSt} as Proposition 7.39 or at \cite{BrFeHo} as Theorem A.2.5.\\
\indent Consequently, \eqref{proof.estimates.1}$_2$ yields
$$
\mathbb E\left[\|\vu\|^r_{L^\infty((0,T);X_n)}\right]\leq C.
$$
We deduce similarly to \cite[Chapter 4.1.2.9]{BrFeHo} that
$$
\mathbb E \left[ \|\vu\|^r_{C^\beta([0,T);X_n)}\right] \leq C
$$
for some $\beta>0$.

Everything is prepared for the use of the Jakubowski-Skorokhod theorem (see \cite[Theorem 2.7.1]{BrFeHo}). We define
\begin{multline*}
\mathcal X = \mathcal X_{\varrho_0}\times \mathcal X_{\vu_0} \times\mathcal X_\varrho\times \mathcal X_\vu\times \mathcal X_W \\ = C^{2+\nu}(\Omega) \times X_n\times C^\iota([0,T],C^{2+\iota}(\Omega))\times C^\kappa([0,T];X_n)\times C([0,T];\mathfrak U_0)
\end{multline*}
for some $\iota,\kappa >0$. Let $\varrho_\mu$ and $\vu_\mu$ be a solution to \eqref{proof.compressible1} and \eqref{proof.momentum1} for some $\mu>0$. The laws of $\{(\varrho_0,\vu_0,\varrho_\mu,\vu_\mu,W)\}_{\mu>0}$ are tight in $\mathcal X$.\\
Indeed, the ball $B_L\subset C^\nu ([0,T]; C^{2+\nu}(\Omega))$ is compact in $\mathcal X_\varrho$ for every $0<\iota<\nu$ and thus
$$
\mathcal L_{\varrho_\mu}(B_L) = \mathbb P(\varrho_\mu \in B_L) = 1
$$
assuming $L$ is sufficiently large due to \eqref{proof.odhad.rho1} (note the constants there are not random variables). Consequently, the laws of $\{\varrho_\mu\}_{\mu>0}$ are tight in $\mathcal X_\varrho$.\\
Similarly, the ball $B_M\subset C^\beta([0,T]; X_n)$ centered at zero with radius $M>0$ is compact in the space $C^\kappa([0,T]; X_n)$ for all $0<\kappa<\beta$. With the help of the Chebyshev inequality we get
\begin{multline*}
\mathcal L_{\vu_\mu}(B_M) = \mathbb P (\vu_\mu \in B_M) = 1-\mathbb P(\|\vu_\mu\|_{C^\beta([0,T];X^n)}>M)\\ = 1-\frac 1{M^r} \mathbb E\left(\|\vu_\mu\|_{C^\beta([0,T];X_n)}^r\right) \geq 1-\frac C{M^r}
\end{multline*} 
and, consequently, the laws of $\{\vu_\mu\}_{\mu>0}$ are tight in $\mathcal X_\vu$.\\
Finally, the law of  $\{W\}$ is tight in $\mathcal X_W$ since it is a singleton. The same reasoning gives that the law of $\varrho_0$ is tight in $C^{2+\nu}(\Omega)$ and the law of $\vu_0$ is tight in $X_n$.\\
According to the Jakubowski-Skorokhod theorem we get the existence of a stochastic basis $(\mathcal O, \mathfrak F, (\mathfrak F_t)_{t\geq 0}, \mathbb P)$ and  random variables $\tilde \varrho_0$, $\tilde \vu_0$, $\tilde \varrho_\mu$, $\tilde \vu_\mu$, $\tilde W_\mu$ which coincide in law to (respectively) $\tilde \varrho_0$, $\tilde \vu_0$, $\varrho_\mu$, $\vu_\mu$, and $\tilde W$. Moreover, there exist a pentad $(\tilde \varrho_0,\tilde \vu_0,\tilde \varrho,\tilde \vu, \tilde W)$ such that 
$$
(\tilde \varrho_0,\tilde \vu_0,\tilde \varrho_\mu, \tilde \vu_\mu, \tilde W_\mu)\to (\tilde \varrho_0,\tilde \vu_0,\tilde \varrho,\tilde \vu, \tilde W)
$$
$\mathbb P$-almost surely in the topology of $\mathcal X$. In particular, $\tilde \varrho$ and $\tilde \vu$ satisfy 
\begin{equation}\label{proof.compressible2}
\pat \tilde \varrho + \diver (\tilde\varrho[\tilde\vu]_R) =\varepsilon \Delta \tilde \vu.
\end{equation}
$\mathbb P$-almost surely with the initial condition $\tilde \varrho(0) = \tilde \varrho$. Next, it satisfies 
\begin{multline}\label{proof.momentum2}
-\int_0^T\pat \Phi\int_\Omega \tilde \varrho \tilde \vu\,  \vphi\, \di x \di t - \Phi(0)\int_\Omega \tilde\varrho_0\tilde\vu_0\, \vphi\, \di x \\=  \int_0^T\Phi\int_\Omega \left(\tilde \varrho [\tilde \vu]_R\otimes \tilde \vu + \chi(\|\tilde \vu - R\|_{X_n}) p(\varrho)\mathbb I\right):\nabla \vphi\, \di x \di t\\ - \int_0^T\Phi \int_\Omega (\mathbb S(D\tilde \vu):\nabla\vphi - \varepsilon \tilde \varrho \tilde \vu \, \Delta \vphi)\, \di x \di t \\
+ \int_0^T\Phi \int_\Omega \tilde\varrho \Pi_n 
[\mathbb F(\tilde \varrho, \tilde \vu)]\cdot \vphi\, \di x \di W
\end{multline}
for each $\Phi \in C_c([0,T))$ and $\vphi\in X_n$ $\mathbb P$-almost surely -- the strong convergence of $\tilde \varrho_\mu$ and $\tilde \vu_\mu$ is sufficient to obtain the convergence of the stochastic process -- we refer to \cite[Lemma 2.6.6 \& proof of Proposition 4.3.15]{BrFeHo}.

Next, we pass with $\varepsilon \to 0$ similarly as above. In particular, let $\varrho_\varepsilon$ and $\vu_\varepsilon$ be the solution to \eqref{proof.compressible2} and \eqref{proof.momentum2} constructed above and emanating from the initial conditions $\varrho_0$ and $\vu_0$. The energy inequality obtained similarly to \eqref{proof.energy2} has the form
\begin{multline}\label{proof.energy3}
\mathbb E \left[\left|\sup_{\tau\in [0,T]} \int_\Omega \frac 12 \varrho_\varepsilon |\vu_\varepsilon|^2 + P(\varrho_\varepsilon)\, \di x \right|^r\right] \\
+ \mathbb E\left[\left|\int_0^T\int_\Omega \mathbb S:\nabla\vu_\varepsilon + \varepsilon \varrho_\varepsilon |\nabla\vu_\varepsilon|^2 + \varepsilon P'' (\varrho_\varepsilon)|\nabla\varrho_\varepsilon|^2\, \di x\di t\right|^r\right]\\
\leq C\, \mathbb E\left[\left|\int_\Omega \frac 12 \varrho_0 |\vu_0|^2 + P(\varrho_0)\, \di x \right|^r + 1\right]
\end{multline}
whenever $r\geq 2$ and one gets also the existence of $C$ such that
\begin{equation}\label{proof.odhad.rho2}
\frac 1C \leq \varrho (t,x) \leq C, \quad \mbox{for all }(t,x)\in [0,T)\times \Omega.
\end{equation}
Consequently, \eqref{proof.energy3} and \eqref{proof.odhad.rho2} yield
\begin{equation*}
\begin{split}
\mathbb E\left[\left|\int_0^T\int_\Omega F(D\vu_\varepsilon)\, \di x \di t\right|^r\right]&\leq C,\\
\mathbb E\left[\|\vu_\varepsilon\|_{L^\infty((0,T); X_n)}^r\right]&\leq C
\end{split}
\end{equation*}
and also 
\begin{equation}\label{proof.estimates.2a}
\mathbb E\left[ \|\vu\|_{C^\beta ([0,T];X_n)}^r\right]\leq c
\end{equation}
for some $\beta>0$.
We consider the sequence (labeled by $\varepsilon>0$)
$$
(\varrho_0,\vu_0,\varrho_\varepsilon,\vu_\varepsilon,W)
$$
in the space 
\begin{multline*}
\mathcal X = \mathcal X_{\varrho_0}\times \mathcal X_{\vu_0}\times\mathcal X_\varrho \times \mathcal X_\vu \times \mathcal X_W \\
=C^{2+\nu}(\Omega)\times X_n\times L^q(((0,T)\times \Omega), w)\\
\times C^\kappa([0,T]; X_n)\times C([0,T]; \mathfrak U_0)
\end{multline*}
for an arbitrary $q\in (1,\infty)$. Again, this sequence is tight:

We recall that $B_L = \{\vu\in L^q((0,T)\times \Omega),\, \|\vu\|_{L^q((0,T)\times\Omega)}\leq L\}$ is compact in the weak topology for every $M>0$ and thus the estimate \eqref{proof.odhad.rho2} yields
$$
\mathcal L_{\varrho_\varepsilon}  = \mathbb P(\varrho_\varepsilon \in B_L) = 1
$$
once $L$ is sufficiently large.\\
Similarly as before, the ball $B_M\subset C^\beta([0,T];X_n)$ is compact in $C^\kappa([0,T];X_n)$ whenever $\kappa<\beta$ and the Chebyshev inequality with \eqref{proof.estimates.2a} give
\begin{multline*}
\mathcal L_{\vu_\varepsilon}(B_M) = \mathbb P(\mu_\varepsilon \in B_M) = 1-\mathbb P(\|\mu_\varepsilon\|_{C^\beta([0,T];X_n)}>M)\\
\geq 1 - \frac 1{M^r}\mathbb E\left(\|\mu_\varepsilon\|_{C^\beta([0,T];X_n)}^r\right)\geq 1 - \frac C{M^r}
\end{multline*}
which implies the desired tightness.

The set $\{\varrho_0,\vu_0,W\}$ is a singleton in $\mathcal X_{\varrho_0}\times \mathcal X_{\vu_0}\times \mathcal X_W$ and therefore it is tight.

As above, the Jakubowski-Skorokhod theorem gives the existence of random variables $\tilde \varrho_0$, $\tilde \vu_0$, $\tilde\varrho_\varepsilon$, $\tilde \vu_\varepsilon$, and $\tilde W$ which coincide in law (respectively) to $\varrho_0$, $\vu_0$, $\varrho_\varepsilon$, $\vu_\varepsilon$, and $W$ such that 
$$
(\tilde \varrho_0,\tilde \vu_0,\tilde \varrho_\varepsilon,\tilde \vu_\varepsilon,\tilde W)\to (\tilde \varrho_0,\tilde \vu_0,\tilde \varrho,\tilde \vu, \tilde W)\ \mbox{ in } \mathcal X
$$
for some pentad $(\tilde \varrho_0,\tilde \vu_0,\tilde \varrho, \tilde \vu, \tilde W)$. The pentad satisfies
\begin{itemize}
\item the continuity equation 
\begin{equation}\label{proof.compressible3}
\pat \tilde \varrho + \diver(\tilde \varrho[\tilde \vu]_R) = 0
\end{equation}
in a weak sense $\mathbb P$-almost surely with initial condition $\tilde\varrho(0) = \tilde\varrho_0$ $\mathbb P$-almost surely. The zero on the right hand side follows from \eqref{proof.energy3} which yields 
$$
\mathbb E \left[\left|\int_0^T\int_\Omega \varepsilon|\nabla \tilde \varrho_\varepsilon |^2\, \di x \di t\right|^r\right]\leq C,
$$
since $P''(\tilde\varrho_\varepsilon)\geq C$ as long as $\tilde \varrho_\varepsilon$ and $\tilde\varrho_\varepsilon^{-1}$ are bounded. Therefore
\begin{equation*}
\mathbb E \left[\left|\int_0^T\int_\Omega \varepsilon \nabla \tilde \varrho_\varepsilon \nabla \varphi\, \di x \di t \right|\right]\leq \sqrt\varepsilon\,  \mathbb E\left[\int_0^T\int_\Omega \varepsilon |\nabla \tilde \varrho _\varepsilon |^2\, \di x \di t\right] \to 0,
\end{equation*}
whenever $\varphi \in C^{\infty}([0,T]\times \Omega)$ is such that $\|\nabla \varphi\|_{L^2((0,T)\times \Omega)} = 1$. Consequently, the right hand side of \eqref{proof.compressible2} tends to zero almost surely.

It is also worthwhile to mention that we can multiply \eqref{proof.compressible2} by $(1+ \log \tilde \varrho_\varepsilon)$ and \eqref{proof.compressible3} by $(1+\log \tilde \varrho)$ to learn
\begin{equation*}
\begin{split}
\pat \left(\tilde \varrho_\varepsilon \log \tilde \varrho_\varepsilon \right) + \diver (\tilde \varrho_\varepsilon \log \tilde \varrho_\varepsilon [\tilde\vu_\varepsilon]_R) \qquad &\\ + \tilde \varrho_\varepsilon\diver \tilde \vu_\varepsilon & = \varepsilon \left(\Delta(\tilde \varrho_\varepsilon \log \tilde \varrho_\varepsilon) - \frac{|\nabla \tilde \varrho_\varepsilon|^2}{\tilde \varrho_\varepsilon}\right)\\
\pat \left(\tilde \varrho \log \tilde \varrho \right) + \diver \left(\tilde \varrho \log \tilde \varrho \tilde \vu \right) + \tilde \varrho \diver \tilde \vu & = 0.
\end{split}
\end{equation*}
Integrating the above equations over $\Omega$ and subtracting them, we deduce
$$
\int_\Omega\overline{\tilde \varrho \log \tilde \varrho}(t,\cdot)\, {\rm d}x \leq \int_\Omega \tilde \varrho \log \tilde \varrho(t,\cdot)\, \di x
$$
for almost all $t\in (0,T)$ where $\overline{f(\tilde \varrho)}$ denotes the weak limit of $f(\tilde \varrho_\varepsilon)$. This is sufficient to claim that $\tilde \varrho_\varepsilon \to \tilde \varrho$ almost everywhere and, consequently, $\tilde \varrho_\varepsilon \to \tilde \varrho$ strongly in $L^p((0,T)\times \Omega)$ for every $p\in [1,\infty)$ almost surely.
\item The momentum equation
\begin{multline*}
-\int_0^T\pat \Phi\int_\Omega \tilde \varrho \tilde \vu\,  \vphi\, \di x \di t - \Phi(0)\int_\Omega \varrho_0\vu_0\, \vphi\, \di x \\=  \int_0^T\Phi\int_\Omega \left(\tilde \varrho [\tilde \vu]_R\otimes \tilde \vu + \chi(\|\tilde \vu - R\|_{X_n}) p(\varrho)\mathbb I\right):\nabla \vphi\, \di x \di t\\ - \int_0^T\Phi \int_\Omega \mathbb S(D\tilde \vu):\nabla\vphi\, \di x \di t 
+ \int_0^T\Phi \int_\Omega \tilde\varrho \Pi_n 
[\mathbb F(\tilde \varrho, \tilde \vu)]\cdot \vphi\, \di x \di W
\end{multline*}
holds for each $\Phi\in C_c([0,T))$ and $\vphi\in X_n$ $\mathbb P$-almost surely.
\item We can neglect certain positive terms on the left hand side of \eqref{proof.energy3} (recall that convex continuous function is weakly sequentially lower semicontinuous) to deduce
\begin{multline}\label{proof.energy4}
\mathbb E \left[\left|\sup_{\tau\in [0,T]} \int_\Omega \frac 12 \varrho |\vu|^2 + P(\varrho)\, \di x \right|^r\right] 
+ \mathbb E\left[\left|\int_0^T\int_\Omega \mathbb S:\nabla\vu\, \di x\di t\right|^r\right]\\
\leq C \mathbb E\left[\left|\int_\Omega \frac 12 \varrho_0 |\vu_0|^2 + P(\varrho_0)\, \di x \right|^r + 1\right]
\end{multline}
\end{itemize}
with $R$-independent constant $C$ on the right hand side. It is worthwhile to mention that the above solution is unique. The stochastic terms in the above system are treated as above -- i.e. the strong convergence of both $\tilde \varrho_\varepsilon$ and $\tilde \vu_\varepsilon$ is sufficient. 

It remains to tend with $R$ to infinity. This is done by a suitable use of the stopping-time argument (see \cite[Chapter 4.2]{BrFeHo}). In particular, let $\varrho_R$ and $\vu_R$ be the unique solution obtained in the previous paragraph and let $\di_t$ denote the material derivative with the velocity field $[\vu]_R$. Then, \eqref{proof.compressible3} yields 
$$
\di_t \varrho_R = -\vr_R \diver[\vu_R]_R
$$
and we deduce 
\begin{equation*}
 \underline \varrho e^{-\int_0^t \|\diver [\vu_R]_R\|_{L^\infty(\Omega)}\, \di t}\leq \varrho(t,x)\leq \overline \varrho e^{\int_0^t\|\diver[\vu_R]_R\|_{L^\infty(\Omega)}\, \di t},
\end{equation*}
where $\underline \varrho \leq \varrho_0(x)\leq \overline \varrho$ for every $x\in \Omega$. Since all norms on $X_n$ are continuous, the Young inequality and \eqref{F.potential.2} yield
\begin{multline*}
\|\diver [\vu_R]_R\|_{L^\infty(\Omega)}\leq \|D \vu_R\|_{L^\infty(\Omega)} \sim \|D \vu_R\|_{L^1(\Omega)}\\ \leq C g^*(1) + \int_\Omega F(D\vu_R)\, \di x.
\end{multline*}
We deduce from  \eqref{proof.energy4} that 
\begin{equation*}
\mathbb E \left[e^{-\int_0^T\int_\Omega F(D\vu_R)\, \di x \di t} \sup_{t\in (0,T)} \|\vu_R\|_{L^2(\Omega)}^2 \right] \leq C.
\end{equation*}
We define three events:
\begin{equation*}
\begin{split}
A& = \left[e^{-\int_0^T\int_\Omega F(D\vu_R)\, \di x \di t} \sup_{t\in (0,T)} \|\vu_R\|_{L^2(\Omega)}^2 \right],\\
B& = \left[\int_0^T\int_\Omega F(D\vu_R)\, \di x \di t \leq b_R\right],\\
S& = \left[\sup_{t\in (0,T)}\|\vu_R\|_{L^2(\Omega)}^2\leq a_Re^{b_R}\right],
\end{split}
\end{equation*}
where $a_R$ and $b_R$ are sequences satisfying $\lim a_R = \lim b_R = \infty$ and $a_r e^{b_R} = R$. Next, because $S\supset A\cap B$, the Chebyshev inequality together with \eqref{proof.energy4} give
\begin{equation}\label{proof.stopping.time}
\mathbb P(S)\geq \mathbb P(A) + \mathbb P(B) - 1 \geq 1-\frac C{a_R} - \frac C{b_R} \to 1.\end{equation}
We define the stopping time 
$$
\tau_R = \inf \{t\in [0,T],\ \|\vu_R\|_{X_n} >R\}.
$$
Recall that $[\vu_R(t,\cdot)]_R = \vu_R(t,\cdot)$ whenever $t<\tau_R$. We have just deduced that $\mathbb P(\sup_{R\in \mathbb N}\tau_R = T) = 1$ due to \eqref{proof.stopping.time}.
\section{Measure-valued martingale solution}\label{m-v.m.s}
First, we construct the laws of initial conditions $\Lambda_n$ in such a way that $\Lambda_n\to \Lambda$ in certain sense and $\Lambda_n$ satisfies the assumptions from Theorem \ref{Galerkin.existence} and $\Lambda$ is a probability measure on $L^1(\Omega)\times L^1(\Omega)$ such that
\begin{equation}\label{lambda.ini.con.1}
\Lambda\{\varrho\geq 0\} = 1,\ \Lambda\left\{0<\underline \varrho \leq \int_\Omega \varrho \, \di x \leq \overline \varrho <\infty\right\} = 1
\end{equation}
and
\begin{equation}\label{lambda.ini.con.2}
\int_{L^1(\Omega)\times L^1(\Omega)}\left|\int_\Omega \left[\frac 12 \frac{(\varrho\vu)^2_0}{\varrho_0} + P(\varrho)\right]\, \di x \right|^r\, \di \Lambda \leq C
\end{equation}
for some $r\geq 4$. In particular, let $\varrho_0$ and $\vu_0$ be the random variable on $L^1(\Omega)\times L^1(\Omega)$ whose law is given by $\Lambda$. Then we construct  $\varrho_{0,n}$ and $\vu_{0,n}$ as the initial conditions satisfying
$$
(\varrho_{0,n},\vu_{0,n})\in C^{2+\nu}(\Omega)\times X_n,\quad \frac 1 {\varrho_{0,n}}\in L^\infty(\Omega)
$$
and the following convergencies hold true:
\begin{equation*}
\begin{split}
\varrho_{0,n}&\wto\varrho_0\ \mbox{in }L^P(\Omega),\\ \varrho_{0,n}\vu_{0,n}&\wto (\varrho\vu)_0\ \mbox{in }L^G(\Omega),\\
\int_\Omega \left[\frac 12 \varrho_{0,n}\vu_{0,n}^2 + P(\varrho_{0,n})\right]\, \di x&\to \int_\Omega \left[\frac 12 \frac{(\varrho\vu)_0^2}{\varrho_0} + P(\varrho_0)\right]\, \di x
\end{split}
\end{equation*}
$\mathbb P$-almost surely where $G$ is some N-function. Consequently, we take $\Lambda_n$ to be the law of $\varrho_{0,n}$ and $\vu_{0,n}$.

Throughout this section we work with the sequence $(\varrho_n, \vu_n, \mathbb S_n)$ of solutions constructed in the previous section and emanating from $\varrho_{0,n}$ and $\vu_{0,n}$ whose laws are given by $\Lambda_n$ and we let $n\to \infty$. In order to reach this goal, we use the advantage of stochastic Young measures (as stated in Lemma 2.2 and Lemma 2.3 in \cite{HoKoSa}). In particular, we use \eqref{no.phi.energy} in order to deduce the following bounds:
\begin{equation}\label{approx.bounds}
\begin{split}
\mathbb E\left(\left|\int_Q F^*(\mathbb S_n)\, {\rm d}x{\rm d}t\right|^r\right)&\leq C,\\
\mathbb E\left(\|P(\varrho_n)\|^r_{L^\infty(L^1(\Omega))}\right)& \leq C,\\
\mathbb E\left(\left|\int_Q F(\mathbb D u_n)\, {\rm d}x{\rm d}t\right|^r\right)& \leq C,\\
\end{split}
\end{equation}
where $C$ depends on initial data, right hand side and $\Omega$.
\begin{Lemma}\label{Lemma.4}
Let $F$ satisfy \eqref{F.potential.2}. Assume $\vu_n\restriction_{\partial\Omega} = 0$ and let \eqref{approx.bounds}$_3$ hold. Then there is an $N$-function $g_0:[0,\infty)\to [0,\infty)$ such that
$$
\mathbb E\left(\left|\int_Q g_0(|\vu_n|) \, \di x \di t\right|^r\right)\leq C.
$$
\end{Lemma}

\begin{proof}
Let $\{\vu_{n,1}\}$ be the sequence of the first coordinates of $\vu_n$. For every $x = (x_1,x_2, x_3)\in \partial \Omega$ let $x_b = (x_{1b},x_2,x_3),\ x_{1b}<x_1$ be the point of the boundary such that all points with coordinates $(s,x_2,x_3),\ s\in (x_{1b}, x_1)$ belong to $\Omega$. Assume without loss of generality that $\Omega$ is small enough that $x_{1b}< x_1<x_{1b} + 1$ (otherwise, there will appear some additional constant in the following computation) and, for the purpose of this proof, assume that $\partial_1 \vu_{n,1} (s,x_2,x_3) = 0$ whenever $(s,x_2,x_3)\notin \Omega$.\\
\indent Let $g_1$ be an (arbitrary) $N$-function. Then the Jensen inequality yield 
\begin{multline}\label{poincare.like.1}
\int_\Omega g_1(|\vu_{n,1}|)\, \di x  \leq \int_\Omega g_1\left( \int_{x_{1b}}^{x_1}| \partial_1 \vu_{n,1}|\, \di s\right)\di x_1\di x_2\di x_3 \\ \leq \int_\Omega g_1\left(\int_{x_{1b}}^{x_{1b}+1}|\partial_1 \vu_{n,1}|\, \di s\right)\, \di x_1\di x_2\di x_3\\ \leq \int_\Omega\left(\int_{x_{1b}}^{x_{1b}+1} g_1(|\partial_1 \vu_{n,1}|)\, \di s\right)\, \di x_2\di x_3  \leq \int_\Omega g_1\left(|\mathbb D \vu_n|\right)\, \di x\\
\leq \int_\Omega F(D\vu_n)\, \di x
\end{multline}
assuming $g_1 = g$.
This computation can be repeated in order to get $g_2$ and $g_3$ for the remaining variables. Furthermore, there exists a constant $c_{2,\infty}\in (0,\infty)$ such that any $N$-function $g_0$ satisfies
\begin{multline*}
\int_\Omega g_0(|\vu_n|)\, \di x \leq \int_\Omega g_0(c_{2,\infty} \max \{|\vu_{n,1}|,|\vu_{n,2}|, |\vu_{n,3}|\})\, \di x \\ 
\leq \int_\Omega g_0(c_{2,\infty}|\vu_{n,1}|)\, \di x + \int_\Omega g_0(c_{2,\infty}|\vu_{n,2}|)\, \di x + \int_\Omega g_0(c_{2,\infty}|\vu_{n,3}|)\, \di x
\end{multline*}
and the choice $g_0(t) = \max \left\{g_1\left(\frac t {c_{2,\infty}}\right), g_2\left(\frac t {c_{2,\infty}}\right), g_3\left(\frac t {c_{2,\infty}}\right)\right\}$ together with \eqref{poincare.like.1} yield
\begin{multline*}
\int_\Omega g_0(|\vu_n|)\, \di x \leq \int_\Omega g_1(|\vu_{n,1}|)\, \di x + \int_\Omega g_2(|\vu_{n,2}|)\, \di x  + \int_\Omega g_3(|\vu_{n,3}|)\, \di x \\
\leq 3 \int_\Omega F(D\vu_n)\, \di x,
\end{multline*}
which completes the proof.
\end{proof}

We state the path space $\mathcal X$ as
\begin{multline*}
\mathcal X = \mathcal X_{\varrho_0}\times \mathcal X_{(\varrho\vu)_0}\times \mathcal X_{\varrho\vu\otimes\vu}\times \mathcal X_{p(\varrho)}\times \mathcal X_E \times \mathcal X_\nu \\
=L^P(\Omega,w)\times L^G(\Omega,w)\times  L^\infty(0,T;\mathcal M(\Omega),w^*)\\   \times   L^\infty(0,T;\mathcal M(\Omega),w^*)  \times L^\infty((0,T)\times\Omega; \mathbb P(\mathbb R^{22}),w^*)
\end{multline*}
and we work with the sequence $(\varrho_{0,n},(\varrho_{0,n}\vu_{0,n}), \varrho_n\vu_n\otimes \vu_n, p(\varrho_n), E_n, \nu_n)\in \mathcal X$ where
$$
\nu_{n} = \delta_{\varrho_n}\times \delta_{\vu_n}\times \delta_{\mathbb S_n}\times \delta_{D\vu_n}.
$$
and $E_n$ is the abbreviation for
$$
E_n = \frac 12 \varrho_n|\vu_n|^2 + P(\varrho_n)
$$
Clearly, the boundedness and convergences of $\varrho_{0,n}$ and $\varrho_n\vu_{0,n}$ implies the tightness of their laws in $\mathcal X_{\varrho_0}$ and $\mathcal X_{(\varrho\vu)_0}$.\\
Next, the sequence $\nu_n$ is tight in $\mathcal X_\nu$. Indeed, this follows as the set 
\begin{multline*}
B_R:= \Big\{\nu\in L^\infty((0,T)\times\Omega; \mathbb P(\mathbb R^{22}),w^*),\,\\ \int_0^T\int_\Omega \int_{\mathbb R^{22}}\left[ |r| + |\vw| + |S| + |D|\right]\, \di \nu \, \di x\di t \leq R\Big\}
\end{multline*}
is compact (see \cite[Proposition 2.8.5]{BrFeHo}). Next, we have due to the Chebyshev and Young inequalities
\begin{multline*}
\mathcal L(\nu_n) (B_R^C) = \mathbb P\left(\int_0^T\int_\Omega \int_{\mathbb R^{22}} \left[|r| + |\vw| + |S| + |D|\right]\, \di \nu_n \, \di x \di t\geq R\right)  \\
\leq \mathbb P\left(\int_0^T \int_\Omega \left( P(\varrho_n) + g(\vu_n) + F^*(\mathbb S_n) + F(D\vu_n)\right)\, \di x \di t  \geq R - C\right)\\
\leq \frac C{R-C}
\end{multline*}
which gives the desired tightness.\\
Finally, the bounded sets in $L^\infty(0,T,\mathcal M(\Omega))$ are relatively compact in weak star topology and therefore the appriori bounds directly yield the desired compactness of the laws of $\varrho_n\vu_n\otimes \vu_n$, $p(\varrho_n)$, and $E_n$.

As a matter of fact, we deduce the existence of the sequence of representations $(\tilde \varrho_0,\tilde{(\varrho\vu)}_0,\widetilde{(\varrho_n\vu_n\otimes\vu_n)},\widetilde{p(\varrho_n)} ,\tilde E_n,  \tilde \nu_n)$ which converges in the appropriate topologies to $(\tilde \varrho_0,\tilde{(\varrho\vu)}_0,\widetilde{(\varrho\vu\otimes\vu)},\widetilde{p(\varrho)} ,\tilde E,  \tilde \nu)$. 
As both $P$ and $F$ are Young functions we use Theorem \ref{thm:hofmanova1} and Lemma \ref{Lemma.4} in order to deduce that the parametrized measure $\nu_{\omega,t,x}$ supported on $[0,\infty)\times \mathbb R^3\times \mathbb R^{3\times 3}_{\sym}\times \mathbb R^{3\times 3}_{\sym}$ satisfies
\begin{equation}\label{mw.abbreviation}
\begin{split}
\tilde\varrho_n(t,x)&\wto \langle \nu_{\omega,t,x}, r\rangle =:\varrho(t,x),\\
\tilde\vu_n(t,x) & \wto \langle \nu_{\omega,t,x}, \vw\rangle =: \vu(t,x),\\
\tilde{\mathbb S}_n(t,x) & \wto \langle \nu_{\omega,t,x}, S\rangle =:\mathbb S(t,x),\\
D \tilde u_n &\wto \langle\nu_{\omega,t,x},D\rangle=:\overline{Du},
\end{split}
\end{equation}
where $(r,\vw,S, D)\in ([0,\infty)\times \mathbb R^3\times \mathbb R^{3\times 3}_{\sym}\times \mathbb R^{3\times 3}_{\sym})$ are dummy variables for the density, velocity and stress tensor and $\tilde \varrho_n$, $\tilde \vu_n$, $\tilde {\mathbb S}_n$ are suitable quantities defined on the standard probability space which coincide with $\varrho_n$, $\vu_n$ and $\mathbb S_n$ in law.

The limit in the energy inequality yields
\begin{multline*}
\mathbb E \left[ \left|
\sup_{\tau \in [0,T]} \int_\Omega \left\langle \nu_{\omega, t,x}, \frac 12 r |\vw|^2 + P(r)\right \rangle\, {\rm d}x + \sup_{\tau\in [0,T]}\int_\Omega \di \mathcal D_t
\right|^r
 \right.\\
\left.+ 
 \left|
 \int_0^T\int_\Omega \left\langle\nu_{\omega,t,x}, F^*(S) + F(D)\right\rangle\, \di x \di t
 \right|^r
\right]\\
\leq C\mathbb E\left[\left|
\int_\Omega \frac 12 \varrho_0 |\vu_0|^2 + P(\varrho_0)\, \di x
\right|^r + 1
\right],
\end{multline*}
where $\mathcal D_t$ is the defect measure arising as
$$
\mathcal D_t  = \lim \left(\frac 12 \varrho_n |\vu_n|^2 + P(\varrho_n)\right) - \left\langle \nu_{\omega, t,x}, \frac 12 r|\vw|^2+ P(r)\right \rangle,
$$
where the limit is understood in the weak sense.

Consequently, 
\begin{equation*}
\begin{split}
\widetilde{\varrho_n \vu_n\otimes \vu_n} &\wto \left\langle \nu, r\vw\otimes \vw\right\rangle + \Theta_t,\\
\widetilde{p(\varrho_n)} &\wto \langle \nu, p(r)\rangle + \Lambda_t
\end{split}
\end{equation*}
for almost all $t\in [0,T]$. The total variations  of measures $\Theta_t$ and $\Lambda_t$ are dominated by the energy inequality -- this is the claim of \cite[Lemma 2.3]{HoKoSa}. In particular, we have 
\begin{equation}\label{tot.var}
\int_\Omega \, \di |\Theta_t| + \int_\Omega \di |\Lambda_t| \leq \int_\Omega \di \mathcal D_t
\end{equation}
\medskip \\

\indent As a matter of fact, the Young measure $\nu$ satisfies the continuity equation
\begin{multline}\label{mw.con}
\int_\Omega\varrho(T)\varphi(T) -\varrho_0\varphi(0)\ {\rm d}x - \int_0^T\int_\Omega \varrho \pat \varphi\, {\rm d}x{\rm d}t\\ - \int_0^T\int_\Omega \langle \nu, r\vw\rangle \nabla \varphi \, {\rm d}x{\rm d}t =0
\end{multline}
for $\varphi\in C^\infty_c(\overline{Q_T})$ $\mathbb P$-a.s. 
\medskip \\
\indent The momentum equation takes the form
\begin{multline}\label{mw.mom}
-\int_0^T\pat \Phi \int_\Omega \langle \nu, r\vw\rangle \cdot \vphi\, {\rm d}x{\rm d}t - \Phi(0)\int_\Omega \varrho_0\vu_0\cdot \vphi(0)\, {\rm d}x \\
= \int_0^T\Phi \int_\Omega \left\langle \nu, r\vw\otimes\vw\right\rangle:\nabla\vphi + \langle \nu,p(r)\rangle\diver\vphi\, {\rm d}x{\rm d}t \\
 + \int_0^T\Phi\int_\Omega \nabla\vphi: \di \Theta_t\di t + \int_0^T\Phi \int_\Omega \diver \vphi\di \Lambda_t\di t\\
- \int_0^T\Phi \int_\Omega\mathbb S:\nabla\vphi \, {\rm d}x{\rm d}t
 + \int_0^T \int_\Omega \Phi \vphi \,{\rm d}x{\rm d}M^1(t)
\end{multline}
where $M^1$ is a martingale satisfying 
\begin{equation}\label{mw.m1mar}
\llangle f(t), M^1 \rrangle  = \sum_{i,j}\left( \sum_{k=1}^\infty \langle \langle \nu_{s,x}^\omega, r {\bf F}_k(r,\vw)\rangle \langle \mathbb D^s f(e_k),g_j\rangle\, {\rm d}s\right) g_i\otimes g_j,
\end{equation}
where $\{g_i\}\subset W^{-r,2}(\Omega)$ is the orthonormal basis for some $r>0$ and $f$ is an arbitrary random process $f$ adapted to $\left(\mathfrak F_t\right)_{t\geq 0}$ satisfying
$$
{\rm d}f = D^d f \, {\rm d}t+\mathbb D^s f\, {\rm d}W,
$$
\begin{multline*}
f\in C([0,T], W^{1,q}(\Omega)\cap C(\Omega)),\ \mathbb E\left[\sup_{t\in [0,T]} \|f\|^2_{W^{1,q}}\right]^q<\infty,\\\ \mathbb P\mbox{-almost surely for all }q\in [1,\infty)
\end{multline*}
and
\begin{multline*}
D^d f\in L^q(\Omega; L^q(0,T;W^{1,q}(\Omega))),\ \mathbb D^s f\in L^2(\Omega;L^2(0,T;L_2(\mathfrak U;L^2(\Omega)))), \\ \left(\sum_{k=1}^\infty |\mathbb D^s f(e_k)|^q\right)^{1/q}\in L^q(\Omega; L^q(0,T;L^q(\Omega))).
\end{multline*}
Indeed, let $M_n$ be the martingale 
$$
M_n = \int_0^t \varrho_n \Pi_n \left[\mathbb F(\varrho_n,\vu_n)\right]\, \di W.
$$
and let $f$ be as above. Then it holds that 
\begin{multline*}
\llangle f,M_n(t) \rrangle \\ = \sum_{i,j=1}^\infty \left(\sum_{k=1}^\infty \int_0^t\langle  \varrho_n \Pi_n\left[{\bf F}_k(\varrho_n,\vu_n)\right], g_i\rangle \langle \mathbb D_t^s f(e_k), g_j\rangle \, \di s\right) g_i\otimes g_j
\end{multline*}
Properties of cross-variation yield that 
$$
f(t) M_n(t) - \llangle f(t),M_n(t)\rrangle
$$
is a martingale and, therefore, 
\begin{multline*}
\mathbb E\Big[L_\tau(\Phi_n) \Big(\langle M_n(t), g_i\rangle \langle f(t),g_j\rangle - \langle M_n(\tau), g_i\rangle \langle f(\tau), g_j\rangle 
\\
 - \sum_{k=1}^\infty \int_\tau^t \langle \varrho_n \Pi_n\left[{\bf F}_k(\varrho_n,\vu_n)\right],g_i\rangle \langle \mathbb D_t^s f(e_k), g_j\rangle \, \di s\Big)\Big ] = 0
\end{multline*}
whenever $0\leq \tau \leq t\leq T$. Here $L_\tau$ is any bounded continuous functional on $\Phi_n = (\varrho_0,\vu_0,\varrho_n\vu_n\otimes \vu_n,p(\varrho_n), E_n,\nu_n, W)$ depending only on the values of $\Phi_n$  restricted to $[0,\tau]$. Due to the convergencies of $\varrho_n$ and $\vu_n$ we get that
$$
L_\tau (\Phi_n)\to L_\tau(\Phi)\ \mathbb P\mbox{-almost surely}.
$$
Recall the assumption on ${\bf F}$ -- especially recall that ${\bf F} = 0$ whenever $\varrho$ is high. Therefore, 
$$
\varrho_n \Pi_n[{\bf F}_k(\varrho_n,\vu_n)] \in L^\infty((0,T)\times \Omega)
$$
uniformly with respect to $n$ and Theorem \ref{thm:hofmanova1} yields that
$$
\varrho_n \Pi_n\left[{\bf F}_k(\varrho_n,\vu_n)\right] \wto \langle \nu, \vG_k(r,\vw)\rangle \ \mbox{ in }L^q((0,T)\times \Omega)
$$
for all $q\in(1,\infty)$. The Sobolev embedding theorem then gives $L^q\subset W^{-m,2}$ for any $m\in \mathbb N$ assuming $q$  is sufficiently large
and thus
$$
\varrho_n \Pi_n\left[{\bf F}_k(\varrho_n,\vu_n)\right] \wto \langle \nu, \vG_k(r,\vw)\rangle \ \mbox{ in }L^2(0,T;W^{-m,2}(\Omega))
$$
for some $m\in \mathbb N$. Consequently,
\begin{multline*}
\mathbb E\Big[L_\tau(\Phi) \Big(\langle M^1(t), g_i\rangle \langle f(t),g_j\rangle - \langle M^1(\tau), g_i\rangle \langle f(\tau), g_j\rangle 
\\
 - \sum_{k=1}^\infty \int_\tau^t \langle \langle\nu,\vG_k(r,\vw)\rangle,g_i\rangle \langle \mathbb D_t^s f(e_k), g_j\rangle \, \di s\Big)\Big ] = 0
\end{multline*}
which gives that 
$$
M^1(t) f(t) - \sum_{i,j=1}^\infty \left(\sum_{k=1}^\infty \left\langle \langle \nu , r {\bf F}_k(r,\vw)\rangle \right \rangle \left\langle \mathbb D^s_t f(e_k),  g_j \right\rangle \, \di s\right) g_i\otimes g_j 
$$
is a martingale. This implies the precise form of the cross-variation for $M^1$ and $f$ and, consequently, also the desired result \eqref{mw.m1mar}.
\medskip \\
\indent In order to deduce the energy inequality, we recall that  $F$ and $F^*$ are convex and therefore sequentially weakly lower semicontinuous. Thus, we deduce from \eqref{no.phi.energy} that
\begin{multline}\label{mw.ene}
\left[\int_\Omega  \left\langle \nu_{\omega,\tau,\cdot},\frac 12r |\vw|^2(\tau,\cdot) + P(r)\right\rangle\, \di x + \int_\Omega \, \di \mathcal D_\tau\right]_{\tau=0}^{\tau=t}\\ + \int_0^\tau \int_\Omega \left\langle \nu, F(S) + F^*( D)\right\rangle\, \di x\di t\\
\leq  \frac 12\int_0^T\int_{\Omega} \sum_{k=1}^\infty \langle \nu_{t,x}^\omega, r^{-1}|{\mathbf G}_k(r,r\vw)|^2\rangle\, {\rm d}x{\rm d}t + \int_0^t \, {\rm d}M^2(t)
\end{multline}
where $\mathcal D_\tau$ is the defect measure.\\
Next, we would like to explain that the limit of $\int_0^t \int_{\Omega} \varrho_n \Pi_n [\mathbb F(\varrho_n,\vu_n)]\cdot \vu_n \,\di x \di W$ which is here denoted by $M^2$ is a martingale. In particular, we show that
$$
\mathbb E[M^2(t)|\mathfrak F_s] = M^2(s)
$$
whenever $0\leq s\leq t\leq T$. As above, we rather show that
$$
\mathbb E\left[L_s(\Phi) (M(t) - M(s))\right] = 0
$$
where $\Phi = (\varrho_0,\vu_0, \varrho \vu\otimes \vu, p(\varrho), E,\nu, W)$ and $L_s$ is any bounded continuous functional that depends only on the values of $\Phi$ restricted to $[0,s]$ on the path space $\mathcal X$. However, we know that
$$
\mathbb E\left[L_s(\Phi_n) \int_s^t \int_{\Omega} \varrho_n \Pi_n [\mathbb F(\varrho_n,\vu_n)]\cdot \vu_n \,\di x \di W\right] = 0
$$
where $\Phi_n = (\varrho_0,\vu_0, \varrho_n \vu_n\otimes \vu_n, p(\varrho_n), E_n,\nu_n, W)$. The convergencies in the appropriate topology of $\mathcal X$ yields that $\Phi_n\to\Phi$ $\mathbb P$-almost surely. This, together with the Vitali convergence theorem allow us to pass to the limit in $n$ and we get the demanded martingale property of the limit $M^2$.
\medskip \\
\indent The above considerations yield the correct notion of a measure-valued martingale solution.
\begin{Definition}
Let $\Lambda$ be a Borel probability measure on $L^\gamma(\Omega)\times L^{2\gamma/(\gamma+1)}(\Omega)$. 
 A trio $[(\pravp,\mathfrak F, (\mathfrak F_t)_{t\geq 0},\mathbb P), \nu_{t,x}^\omega, W]$, where $(\pravp,\mathfrak F, (\mathfrak F_t)_{t\geq 0},\mathbb P)$ is a stochastic basis with a complete right-continuous filtration, $W$ is $(\mathfrak F_t)$-cylindrical Wiener process, and $\nu_{t,x}^w$ is a family of probabilistic measures on $([0,\infty), \mathbb R^3,\mathbb R^{3\times 3},\mathbb R^{3\times 3})$ is called {\emph measure-valued solution } to \eqref{main.sys} if\footnote{We use the abbreviation stated in \eqref{mw.abbreviation}.} the following hold
 \begin{itemize}
 \item the density $\varrho$ satisfies
 $$
 \mathbb E \left[\left(\esssup_{t\in (0,T)} \int_\Omega P(\varrho)\, \di x\right)^r \right] <\infty
 $$
 for all $r\in (1,\infty)$;
 \item there exists an $N$-function $g$ such that
 $$
 \mathbb E\left[\left(\int_0^T\int_\Omega g(|\vu|)\, \di x\di t\right)^r\right]<\infty
 $$
 for all $r\in (1,\infty)$;
 \item the shear-rate $D\vu$ and the stress tensor $\mathbb S$ satisfy
 $$
 \mathbb E \left[\left(\int_0^T\int_\Omega F(D\vu)\, \di x\di t\right)^r + \left(\int_0^T\int_\Omega F^*(\mathbb S)\, \di x \di t\right)^r\right]<\infty
 $$
 for all $r\in (1,\infty)$;
 \item the continuity equation \eqref{mw.con} is satisfied for all $\varphi \in C_c^\infty([0,T)\times \overline\Omega)$ $\mathbb P$-almost surely;
 \item there is a $W^{-m,2}(\Omega)$-valued square integrable continous martingale $M_1$ fulfilling \eqref{mw.m1mar} and measures $\Theta_t$ and $\Lambda_t$ such that the momentum equation \eqref{mw.mom} is satisfied $\mathbb P$-almost surely and for all $\varphi \in C^\infty_c ([0,T]\times \Omega)$;
 \item The total variations of the measures $\Theta_t$ and $\Lambda_t$ satisfy the estimate \eqref{tot.var}  together with the energy inequality \eqref{mw.ene} for almost all $t\in (0,T)$ $\mathbb P$-almost surely.
 \end{itemize}
\end{Definition}
The existence of such solution was shown above. However, our findings are clearly summarized in the following theorem.
\begin{Theorem}
Let the domain $\Omega$ be of class $C^{2+\nu}$ for some $\nu>0$, the right hand side $G$ satisfy \eqref{G.ass.1}, and let the law of the initial data $\Lambda$ satisfies \eqref{lambda.ini.con.1} and \eqref{lambda.ini.con.2}. Next, let $\mathbb S = \mathbb S(D\vu)$ be given as a subdifferential of a convex function $F$ satisfying \eqref{F.potential.2}.

Then there exists a measure-valued martingale solution to \eqref{main.sys} as stated in the above definition.
\end{Theorem}

\bibliographystyle{plain}
\bibliography{literatura}

\begin{thebibliography}{10}

\bibitem{AbFeNo}
Anna Abbatiello, Eduard Feireisl, and Anton\'{\i}n Novotn\'{y}.
\newblock Generalized solutions to models of compressible viscous fluids.
\newblock {\em Discrete Contin. Dyn. Syst.}, 41(1):1--28, 2021.

\bibitem{AdFo}
Robert~A. Adams and John J.~F. Fournier.
\newblock {\em Sobolev spaces}, volume 140 of {\em Pure and Applied Mathematics (Amsterdam)}.
\newblock Elsevier/Academic Press, Amsterdam, second edition, 2003.

\bibitem{basaric}
Danica Basari{\'c}.
\newblock Existence of dissipative (and weak) solutions for models of general compressible viscous fluids with linear pressure.
\newblock {\em J. Math. Fluid Mech.}, 24(2):22, 2022.
\newblock Id/No 56.

\bibitem{Bo}
V.~I. Bogachev.
\newblock {\em Measure theory. {V}ol. {I}, {II}}.
\newblock Springer-Verlag, Berlin, 2007.

\bibitem{BrFeHo}
Dominic Breit, Eduard Feireisl, and Martina Hofmanov\'{a}.
\newblock {\em Stochastically forced compressible fluid flows}, volume~3 of {\em De Gruyter Series in Applied and Numerical Mathematics}.
\newblock De Gruyter, Berlin, 2018.

\bibitem{BrHo}
Dominic Breit and Martina Hofmanov\'{a}.
\newblock Stochastic {N}avier-{S}tokes equations for compressible fluids.
\newblock {\em Indiana Univ. Math. J.}, 65(4):1183--1250, 2016.

\bibitem{DeVo}
A.~Debussche and J.~Vovelle.
\newblock Scalar conservation laws with stochastic forcing.
\newblock {\em J. Funct. Anal.}, 259(4):1014--1042, 2010.

\bibitem{FeLiMa}
E.~Feireisl, X.~Liao, and J.~M\'{a}lek.
\newblock Global weak solutions to a class of non-{N}ewtonian compressible fluids.
\newblock {\em Math. Methods Appl. Sci.}, 38(16):3482--3494, 2015.

\bibitem{FeGwSwWi}
Eduard Feireisl, Piotr Gwiazda, Agnieszka \'{S}wierczewska Gwiazda, and Emil Wiedemann.
\newblock Dissipative measure-valued solutions to the compressible {N}avier-{S}tokes system.
\newblock {\em Calc. Var. Partial Differential Equations}, 55(6):Art. 141, 20, 2016.

\bibitem{FeNoPe}
Eduard Feireisl, Anton\'{\i}n Novotn\'{y}, and Hana Petzeltov\'{a}.
\newblock On the existence of globally defined weak solutions to the {N}avier-{S}tokes equations.
\newblock {\em J. Math. Fluid Mech.}, 3(4):358--392, 2001.

\bibitem{FlGo}
Liviu~C. Florescu and Christiane Godet-Thobie.
\newblock {\em Young measures and compactness in measure spaces}.
\newblock De Gruyter, Berlin, 2012.

\bibitem{HoKoSa}
Martina Hofmanov{\'a}, Ujjwal Koley, and Utsab Sarkar.
\newblock Measure-valued solutions to the stochastic compressible {Euler} equations and incompressible limits.
\newblock {\em Commun. Partial Differ. Equations}, 47(9):1907--1943, 2022.

\bibitem{HoZhZh}
Martina Hofmanov{\'a}, Rongchan Zhu, and Xiangchan Zhu.
\newblock On ill- and well-posedness of dissipative martingale solutions to stochastic 3d {Euler} equations.
\newblock {\em Commun. Pure Appl. Math.}, 75(11):2446--2510, 2022.

\bibitem{Ja}
A.~Jakubowski.
\newblock The almost sure {S}korokhod representation for subsequences in nonmetric spaces.
\newblock {\em Teor. Veroyatnost. i Primenen.}, 42(1):209--216, 1997.

\bibitem{lions}
Pierre-Louis Lions.
\newblock {\em Mathematical topics in fluid mechanics. {V}ol. 2}, volume~10 of {\em Oxford Lecture Series in Mathematics and its Applications}.
\newblock The Clarendon Press, Oxford University Press, New York, 1998.
\newblock Compressible models, Oxford Science Publications.

\bibitem{MNRR}
J.~M\'alek, J.~Ne\v{c}as, M.~Rokyta, and M.~R\r{u}\v{z}i\v{c}ka.
\newblock {\em Weak and measure-valued solutions to evolutionary {PDE}s}, volume~13 of {\em Applied Mathematics and Mathematical Computation}.
\newblock Chapman \& Hall, London, 1996.

\bibitem{MaRa}
J.~M\'{a}lek and K.~R. Rajagopal.
\newblock Compressible generalized {N}ewtonian fluids.
\newblock {\em Z. Angew. Math. Phys.}, 61(6):1097--1110, 2010.

\bibitem{mamontov1}
A.~E. Mamontov.
\newblock On the global solvability of the multidimensional {N}avier-{S}tokes equations of a nonlinearly viscous fluid. {I}.
\newblock {\em Sibirsk. Mat. Zh.}, 40(2):408--420, iii, 1999.

\bibitem{mamontov2}
A.~E. Mamontov.
\newblock On the global solvability of the multidimensional {N}avier-{S}tokes equations of a nonlinearly viscous fluid. {II}.
\newblock {\em Sibirsk. Mat. Zh.}, 40(3):635--649, iii, 1999.

\bibitem{NoSt}
A.~Novotn\'{y} and I.~Stra\v{s}kraba.
\newblock {\em Introduction to the mathematical theory of compressible flow}, volume~27 of {\em Oxford Lecture Series in Mathematics and its Applications}.
\newblock Oxford University Press, Oxford, 2004.

\bibitem{Pedregal}
Pablo Pedregal.
\newblock {\em Parametrized measures and variational principles}, volume~30 of {\em Progress in Nonlinear Differential Equations and their Applications}.
\newblock Birkh\"{a}user Verlag, Basel, 1997.

\bibitem{zhpa}
V.~V. Zhikov and S.~E. Pastukhova.
\newblock On the solvability of the {N}avier-{S}tokes system for a compressible non-{N}ewtonian fluid.
\newblock {\em Dokl. Akad. Nauk}, 427(3):303--307, 2009.

\end{thebibliography}

\end{document}